\documentclass{amsart}
\usepackage{amsmath,amssymb,amsfonts}
\usepackage{amsxtra}
\usepackage{mathtools}
\usepackage{stmaryrd}
\usepackage[all]{xy}

\DeclareMathOperator{\trace}{tr}

\DeclareMathOperator{\Aut}{Aut}
\DeclareMathOperator{\End}{End}
\DeclareMathOperator{\id}{id}

\DeclareMathOperator{\Hom}{Hom}
\DeclareMathOperator{\Pic}{Pic}
\DeclareMathOperator{\Stab}{Stab}
\DeclareMathOperator{\diag}{diag}
\DeclareMathOperator{\Tors}{T}
\DeclareMathOperator{\Fix}{Fix}
\DeclareMathOperator{\str}{str}
\DeclareMathOperator{\diff}{d}
\DeclareMathOperator{\Sym}{S}

\newcommand{\ie}{{\it i.e. }}

\newcommand{\ii}{{\rm i}}
\newcommand{\dual}[1]{{#1}\spcheck}
\newcommand{\abs}[1]{|{#1}|}
\newcommand{\Kummer}[2]{K_{#1}(#2)}

\newcommand{\GAut}{\Aut_{\IZ}} 

\newcommand{\IC}{\mathbb{C}}
\newcommand{\IP}{\mathbb{P}}
\newcommand{\IR}{\mathbb{R}}

\newcommand{\IZ}{\mathbb{Z}}

\newcommand{\cL}{\mathcal{L}}
\newcommand{\cO}{\mathcal{O}}

\newcommand{\kS}{\mathfrak{S}}

\newcommand{\lra}{\longrightarrow}

\newcommand{\genKm}{K_n(A)}
\newcommand{\hilbA}{A^{[n]}}
\newcommand{\quotA}{A^{(n)}}
\newcommand{\fiberA}{A^{((n))}}

\theoremstyle{plain}
\newtheorem{theorem}{Theorem}[section]
\newtheorem{lemma}[theorem]{Lemma}
\newtheorem{proposition}[theorem]{Proposition}
\newtheorem{corollary}[theorem]{Corollary}
\theoremstyle{definition}
\newtheorem{definition}{Definition}[section]
\theoremstyle{remark}
\newtheorem{remark}{Remark}[section]

\begin{document}

\title[Enriques varieties]{Higher dimensional Enriques varieties and\\automorphisms of generalized Kummer varieties}

\author{Samuel Boissi\`ere}

\address{Samuel Boissi\`ere, Laboratoire J.A.Dieudonn\'e UMR CNRS 6621,
         Universit\'e de Nice Sophia-Antipolis, Parc Valrose, F-06108 Nice}

\email{samuel.boissiere@unice.fr}

\urladdr{http://math.unice.fr/$\sim$sb/}

\author{Marc Nieper-Wi{\ss}kirchen}

\address{Marc Nieper-Wi{\ss}kirchen, Lehrstuhl f\"ur Algebra und Zahlentheorie, Universit\"ats\-stra{\ss}e~14, D-86159 Augsburg}
			
\email{marc.nieper-wisskirchen@math.uni-augsburg.de}

\urladdr{http://www.math.uni-augsburg.de/alg/}

\author{Alessandra Sarti}

\address{Alessandra Sarti, Laboratoire de Math\'ematiques et Applications, UMR CNRS 6086,
			Universit\'e de Poitiers, T\'el\'eport 2, Boulevard Marie et Pierre Curie,
			F-86962 Futuroscope Chasseneuil}
			
\email{sarti@math.univ-poitiers.fr}

\urladdr{http://www-math.sp2mi.univ-poitiers.fr/$\sim$sarti/}

\date{\today}



\begin{abstract} 
We define Enriques varieties as a higher dimensional generalization of Enriques surfaces and construct examples by using fixed point free automorphisms on generalized Kummer varieties. We also classify all automorphisms of generalized Kummer varieties that come from an automorphism of the underlying abelian surface.
\end{abstract}

\maketitle

\bigskip

\section{Introduction}

Very much is known on automorphisms of K3 surfaces, but the study of the automorphisms of irreducible holomorphic symplectic manifolds, which are the higher-dimensional analogous of K3 surfaces, is quite recent. Some results are given by Beauville~\cite{B,B3}, Boissi\`ere~\cite{Boissiere}, Boissi\`ere--Sarti~\cite{BS}, Debarre~\cite{Debarre}, Huybrechts~\cite{H} and Oguiso \cite{Oguiso}. In this paper we study generalized Kummer varieties, showing that an automorphism leaving invariant the exceptional divisor is a \emph{natural} automorphism, \ie is induced by an automorphism of the underlying abelian surface. A proof of the corresponding result for the Hilbert schemes of points on K3 surfaces is given in \cite{BS}.

We study also automorphisms of finite order without fixed point on irreducible holomorphic symplectic manifolds. On K3 surfaces, the automorphisms of finite order without fixed points are non-symplectic involutions; the quotients are exactly the Enriques surfaces. This motivates us to introduce in Definition \ref{def:Enriques} the notion of an Enriques variety to generalize the properties of Enriques surfaces in higher dimension. We give in Proposition~\ref{prop:classification} a classification result concerning Enriques varieties and we show in Proposition~\ref{prop:construction} the existence of Enriques varieties of dimension $4$ and $6$ by giving explicit examples as quotients of generalized Kummer varieties by fixed point free natural automorphisms of order $3$ and $4$. This gives a positive answer to a question of Arnaud Beauville asked during the conference \emph{Moduli} in Berlin in 2009. These examples are also constructed by Oguiso--Schr\"oer~\cite{OS} in an independent paper. We thank Arnaud Beauville and Igor Dolgachev for helpful comments and their interest in this work.

\section{Higher dimensional Enriques varieties}

\subsection{Irreducible holomorphic symplectic manifolds}
\label{ss:IHS}

A complex, compact, K\"ahler manifold $X$ is called \emph{irreducible symplectic} if $X$ is simply connected and $H^0(X,\Omega^2_X)$ is spanned by an everywhere non-degenerate two-form, denoted $\omega_X$. In this case, the second cohomology group $H^2(X,\IZ)$ has no torsion and is equipped with a non-degenerate bilinear symmetric form called the Beauville-Bogomolov form~\cite{B3}. We denote by $\cO(H^2(X,\IZ))$ the group of isometries of the lattice $H^2(X,\IZ)$ with respect to this quadratic form. There is a natural map $\Aut(X)\to\cO(H^2(X,\IZ))$, $f\mapsto f^*$; in contrary to the case of K3 surfaces, for $\dim X>2$ this map is not injective in general (though it is indeed injective for Hilbert schemes of points on K3 surfaces \cite{B}), but its kernel is a finite group \cite{H}: an automorphism of $X$ acting trivially on the second cohomology leaves invariant a K\"ahler class and the associated Calabi--Yau metric, so it is an isometry of $X$. Since the group $\Aut(X)$ is discrete and the group of isometries is compact, it is a finite group.

\subsection{Enriques varieties}

Let $S$ be a K3 surface admitting a fixed point free involution $\iota\colon S\to S$. This involution is necessarily non-symplectic and $S$ is projective. The quotient $Y:=S/\langle\iota\rangle$ is an \emph{Enriques surface}, characterized by the numerical conditions $2K_Y=0$ and $\chi(Y,\cO_Y)=1$. It is well known that every Enriques surface can be obtained as such a quotient. They are easy to construct, for example using complete intersections of quadrics in $\IP_5$. One first possible generalization in higher dimension is by using Calabi--Yau manifolds. In odd dimension, the quotient varieties $Y$ obtained as quotient of Calabi--Yau varieties by fixed point free automorphisms satisfy $\chi(Y,\cO_Y)=0$ so would rather correspond to a generalization of bielliptic surfaces. In even dimension, only involutions can act without fixed point on an even dimensional Calabi--Yau manifold since its holomorphic Euler characteristic is two, so the canonical divisor of the quotient can not have order higher that two. Examples are easy to construct: complete intersections of quadrics (see \S\ref{ss:counterexample}) or --- as was pointed out to us by Igor Dolgachev --- generalization in higher dimension of Reye congruences \cite{Cossec} (see \cite{CD} for other constructions).
The problem is to construct quotients by fixed point free higher order automorphisms, in such a way that the canonical divisor has higher order.

Let $X$ be an irreducible holomorphic symplectic manifold of dimension $2n-2$ with $n\geq 2$, and $f$ an automorphism of $X$ of order $d\geq 2$ such that the cycle group $\langle f \rangle$ generated by $f$ acts freely on $X$. Observe that this is possible only if $f$ is purely non-symplectic, \ie there exists a primitive $d$-th root of the unity $\xi$ such that the action of $f$ on the symplectic form $\omega_X$ is given by $f^*\omega_X=\xi\omega_X$. Indeed, otherwise there would exist some integer $i$, $1\leq i\leq d-1$, such that $f^i$ is symplectic (\ie $(f^i)^*\omega_X=\omega_X$), but this would imply that $f^i$ has fixed points: since $H^0(X,\Omega_X^k)$ is zero for odd $k$ and generated by $\omega_X^{k/2}$ for even $k$, the holomorphic Lefschetz number of $f^i$ is:
$$
L(f^i)=\sum_{j=0}^{2n-2} (-1)^j\trace\left((f^i)^*_{|H^0(X,\Omega_X^j)}\right)=\frac{\dim X}{2}+1=n
$$
and by the holomorphic Lefschetz fixed point formula, $L(f^i)\neq 0$ implies that the fixed locus of $f^i$ is non empty. So the group $\langle f \rangle$ acts purely non-symplectically on $X$, and such a group can exist only when $X$ is projective \cite[Proposition 6]{B}. Then $Y{\coloneqq}X/\langle f\rangle$ is smooth and projective, and since the quotient map $\pi\colon X\to Y$ is a non ramified covering of order $d$, the canonical divisor $K_Y$ has order $d$ in $\Pic(Y)$ and $d\cdot\chi(Y,\cO_Y)=\chi(X,\cO_X)=n$ so $d$ divides $n$. Note that the same argument shows that a fixed point free involution on an even-dimensional Calabi--Yau manifold can not act trivially on the holomorphic volume form. This motivates the following definition:

\begin{definition} \label{def:Enriques}
Let $Y$ be a connected, compact, smooth, k\"ahler, complex manifold.
$Y$ is called an \emph{Enriques variety} if there exists an integer $d\geq 2$, called the \emph{index} of~$Y$, such that the canonical
bundle $K_Y$ has order $d$ in the Picard group $\Pic(Y)$ of~$Y$, the holomorphic Euler characteristic
of $Y$ is $\chi(Y,\cO_Y)=1$ and the fundamental group $\pi(Y)$ is cyclic of order $d$.
\end{definition}

Observe that if $Y_1, Y_2$ are Enriques varieties of indices $d_1,d_2$ prime to each other, then $Y_1\times Y_2$ is again an Enriques variety, of index $d_1d_2$. We thus introduce the following definition:

\begin{definition}
An Enriques variety is called \emph{irreducible} if the holonomy representation of its universal cover is irreducible.
\end{definition}

\begin{proposition}\label{prop:classification} \text{}
\begin{enumerate}
\item\label{prop:itemeven} Every Enriques variety is even dimensional.
\item\label{prop:itemCY} Every irreducible Enriques variety of index two is the quotient of an even dimensional Calabi--Yau variety by a fixed point free involution.
\item\label{prop:itemIHS}  Every irreducible Enriques variety of index strictly greater than two is the quotient of an irreducible symplectic holomorphic manifold by an automorphism acting freely.
\item\label{prop:itemproj} Every irreducible Enriques variety is projective.
\item\label{prop:itemirred} Every Enriques variety of prime index is irreducible.
\item\label{prop:itemirredodd} Every Enriques variety of odd index $d$ and dimension $2d-2$ is irreducible.
\end{enumerate}
\end{proposition}

\begin{proof} Let $Y$ be an Enriques variety. The case $\dim Y=2$ is clear so we assume that $\dim Y>2$.
Since $d K_Y=0$, $Y$ admits an order $d$ non ramified covering $X\to Y$ such that $K_X=0$. Since $X$ is simply connected, by the decomposition theorem of Bogomolov it is isomorphic to a product of Calabi--Yau varieties and irreducible symplectic holomorphic varieties. Since the holomorphic Euler characteristic of an odd dimensional Calabi--Yau variety is zero, no such variety appears in the decomposition; this proves (\ref{prop:itemeven}). If moreover $Y$ is irreducible, only one factor occurs and since even dimensional Calabi--Yau varieties have holomorphic Euler characteristic equal to two, only involutions can act without fixed point on them; this proves (\ref{prop:itemCY}), (\ref{prop:itemIHS}) and (\ref{prop:itemproj}). For~(\ref{prop:itemirred}), since the Euler characteristic of $X$ is the product of those of its factors and equals the index of $Y$, only one factors occurs when this index is prime so $Y$ is irreducible. For~(\ref{prop:itemirredodd}), decomposing $X\cong W_1\times\cdots\times W_k$ were the $W_j$ are irreducible holomorphic symplectic varieties of dimension $r_j>2$ , one has $\dim X=2d-2=r_1+\cdots+r_k$ and $\chi(X,\cO_X)=d=\left(\frac{r_1}{2}+1\right)\cdots \left(\frac{r_k}{2}+1\right)$ so if $k\geq 2$, $2d=r_1+\cdots+r_k+\frac{1}{2} r_1r_2+\text{positive terms}$ with $r_1r_2>4$, that is not possible so $k=1$.
\end{proof}

Note that the assertion (\ref{prop:itemirredodd}) would be false for an even index different from two: a counter-example is constructed in \S\ref{ss:counterexample}. In \ref{prop:construction} we construct irreducible Enriques varieties of index three and four, thus proving:

\begin{theorem}
There exists irreducible Enriques varieties of index $2$, $3$ and $4$.
\end{theorem}

It is easy to construct examples of a weaker notion of Enriques varieties by asking only that $\chi(Y,\cO_Y)\neq 0$. Such varieties can appear as intermediate quotients between an irreducible holomorphic symplectic manifold and an Enriques variety. We will call them \emph{weak Enriques varieties}.

\section{Generalized Kummer varieties}

Let $A$ be a complex torus of dimension two and $n\geq 2$. Denote by $\quotA{\coloneqq}~A^n/\kS_n$ the symmetric quotient of $A$, where $\kS_n$ acts by permutation of the factors. Let $\pi\colon A^n\to~\quotA$ be the quotient map, $\Delta$ the large diagonal in $A^n$ and $D{\coloneqq}\pi(\Delta)$ its image in $\quotA$. Let $\hilbA$ be the \emph{Douady space} (or \emph{Hilbert scheme} when $A$ is algebraic) parametrizing zero-dimensional length $n$ analytic subspaces of $A$ and ${\rho\colon\hilbA\to\quotA}$ the \emph{Douady--Barlet} (\emph{Hilbert--Chow} in the algebraic case) morphism. It is a resolution of singularities  and the exceptional divisor $E{\coloneqq}\rho^{-1}(D)$ is irreducible. Let $s\colon\quotA\to A$ be the summation map and consider the composed morphism $S\colon\hilbA\xrightarrow{\rho}\quotA\xrightarrow{s} A$. Set $\genKm{\coloneqq}S^{-1}(0)$. Then $\genKm$ is a smooth, complex, compact, irreducible symplectic holomorphic manifold of dimension $2n-2$. For ${n=2}$, $\Kummer 2 A$ is the Kummer surface associated to $A$. Set $\fiberA{\coloneqq}s^{-1}(0)$. The restriction of the Hilbert--Chow morphism to $\genKm$ is again denoted by ${\rho\colon\genKm\to\fiberA}$. It is a resolution of singularities and the exceptional divisor $E_0{\coloneqq}E\cap\genKm$ is irreducible for $n\geq 3$. We set $D_0{\coloneqq}D\cap\genKm$, $\tilde{s}\colon A^n\to A$ the summation map, $A^n_0{\coloneqq}\tilde{s}^{-1}(0)$, $\Delta_0{\coloneqq}\Delta\cap A^n_0$ so that, for the quotient map $\pi\colon A^n_0\to A^n_0/\kS_n\cong\fiberA$, one has $D_0=\pi(\Delta_0)$.

\subsection{Characterization of the natural automorphisms}

Any biholomorphic map $\psi\colon A\to A$ induces in a natural way an automorphism $\psi^{[n]}\colon\hilbA\to\hilbA$, called \emph{natural} \cite{Boissiere} and the map $\Aut(A)\to\Aut(\hilbA)$, $\psi\mapsto \psi^{[n]}$ is injective. Such an automorphism restricts to an automorphism of $\genKm$ if and only if for any $\xi\in\hilbA$ such that $S(\xi)=0$ one has $S(\psi^{[n]}(\xi))=0$. Recall \cite[\S1.2]{BL} that any biholomorphic map $\psi\colon A\to A$ decomposes in a unique way as $\psi=t_{\psi(0)}\circ h$ where $t_{\psi(0)}$ is the translation by $\psi(0)$ in $A$ and $h\in\GAut(A)$ is a group automorphism, so that $\Aut(A)\cong A\rtimes\GAut(A)$. One can easily see that for any $h\in\GAut(A)$, $h^{[n]}$ restricts to an automorphism of $\genKm$, and for $a\in A$, the translation $t_a$ by $a$ in $A$ induces an automorphism of $\hilbA$ that restricts to $\genKm$ if and only if $a$ is an $n$-torsion point of $A$. Denoting by $\Tors_n(A)$ the $n$-torsion subgroup of $A$, we thus get a well defined morphism $\Tors_n(A)\rtimes\GAut(A)\to\Aut(\genKm)$, $\psi\mapsto \psi^{{\llbracket}n\rrbracket}$, whose image is called the group of \emph{natural automorphisms} of $\genKm$. 

\begin{theorem}\label{th:naturalKummer}
Let $A$ be a complex torus of dimension two and $n\geq 3$. An automorphism of $\genKm$ is natural if and only if it leaves the exceptional divisor globally invariant.
\end{theorem}

\begin{proof}
The condition is clearly necessary; let us show that it is sufficient. Let $f\colon\genKm\xrightarrow{\sim}\genKm$ be an automorphism such that $f(E_0)=E_0$.

\par\noindent{\it Step 1. Universal cover.}\\
The automorphism $f$ induces an automorphism of $\genKm\setminus E_0\cong \fiberA\setminus D_0$. Let $A=\IC^2/\Gamma$ where $\Gamma\subset\IC^2$ is a rank-4 lattice and set
\begin{align*}
A^n_0&{\coloneqq}\left\{(a_1,\ldots,a_n)\in A^n\,|\,\sum_{i=1}^n a_i=0\right\},\\
\Gamma^n_0&{\coloneqq}\left\{(\gamma_1,\ldots,\gamma_n)\in \Gamma^n\,|\,\sum_{i=1}^n \gamma_i=0\right\},\\
(\IC^2)^n_0&{\coloneqq}\left\{(z_1,\ldots,z_n)\in (\IC^2)^n\,|\,\sum_{i=1}^n z_i=0\right\},\\
\end{align*}
with the quotient maps
$$
(\IC^2)^n_0\overset{p}{\lra}(\IC^2)^n_0/\Gamma^n_0\cong A^n_0\overset{\pi}{\lra}A^n_0/\kS_n=\fiberA
$$
and finally set $\Sigma{\coloneqq}p^{-1}(\Delta_0)$. The restricted map 
$$
(\IC^2)^n_0\setminus \Sigma\overset{p}{\lra} A^n_0\setminus \Delta_0\overset{\pi}{\lra}\fiberA\setminus D_0
$$
is a finite unramified Galois covering with Galois group $G{\coloneqq}\Gamma^n_0\rtimes\kS_n$. Precisely, the group
$\kS_n$ acts on $\Gamma^n_0$ (and equally on $(\IC^2)^n_0$) as follows: 
for any $\sigma\in\kS_n$ and $\gamma=(\gamma_1,\ldots,\gamma_n)\in\Gamma^n_0$,
$$
\sigma\gamma{\coloneqq}(\gamma_{\sigma^{-1}(1)},\ldots,\gamma_{\sigma^{-1}(n)}),
$$
the semi-direct product is $(\gamma,\sigma)(\lambda,\tau)=(\gamma+\sigma\lambda,\sigma\tau)$ and $(\gamma,\sigma)\in G$ acts on ${z\in(\IC^2)^n_0}$ by 
$(\gamma,\sigma) z{\coloneqq}\sigma z+\gamma$. Since $\Sigma$ has complex codimension two in $(\IC^2)^n_0$, $(\IC^2)^n_0\setminus\Sigma$ is simply connected, so $\pi\circ p_{|(\IC^2)^n_0\setminus\Sigma}$ is the universal covering of $\fiberA\setminus D_0$. By its universal property, there exists a unique biholomorphic map $$F\colon(\IC^2)^n_0\setminus\Sigma\to(\IC^2)^n_0\setminus\Sigma
$$ 
making the diagram
$$
\xymatrix{(\IC^2)^n_0\setminus\Sigma \ar[r]^-{F}_-{\sim} \ar[d]_p & (\IC^2)^n_0\setminus\Sigma\ar[d]^p\\
A^n_0\setminus\Delta_0\ar[d]_\pi& A^n_0\setminus\Delta_0\ar[d]^\pi\\\fiberA\setminus D_0\ar[r]^-{f}_-{\sim}&\fiberA\setminus D_0}
$$
commutative and with the property that there exists an automorphism $f_*\colon G\to G$ ($G$ is the deck transformation group of the universal covering) such that the map $F$ is $f_*$-equivariant in the following sense:
$$
\forall (\gamma,\sigma)\in G, \forall z\in(\IC^2)^n_0\setminus\Sigma, \quad F((\gamma,\sigma) z)=f_*(\gamma,\sigma) F(z).
$$
Since $\Sigma$ has codimension $2$ in $(\IC^2)^n_0$, by Hartog's theorem $F$ can be extended to a biholomorphic map also denoted $F\colon(\IC^2)^n_0\to(\IC^2)^n_0$ and, by the identity theorem, the $f_*$-equivariance extends to $(\IC^2)^n_0$.

\par\noindent{\it Step 2. Reduction to the case $F(0)=0$.}\\ 
Set $F(0)=(w_1,\ldots,w_{n-1},w_n)$ with $w_i\in\IC^2$ and $\sum_{i=1}^n w_i=0$. Since $F(0)$ is defined up to the action of an element of $G$, by applying a suitable translation by an element $\gamma\in\Gamma^n_0$ one can assume that $(w_1,\ldots,w_{n-1})$ lies in a fundamental domain of $\Gamma$ (this modification of $F$ makes it equivariant for the conjugate $(\gamma,\id)\cdot~f_*(-)\cdot~(\gamma,\id)^{-1}$ of $f_*$, that we continue to write $f_*$ for simplicity). 

Let $\Stab(0)$ be the stabilizer of $0\in(\IC^2)^n_0$ in $G$; then $\Stab(F(0))=f_*(\Stab(0))$. One can easily see that $\Stab(0)=\{(0,\sigma)\in G\,|\, \sigma\in\kS_n\}\cong\kS_n$, so the stabilizer of $F(0)$ is a subgroup of $G$ isomorphic to $\kS_n$. We need the following lemma, whose proof is elementary:

\begin{lemma}\label{lem:stab} Let $H\subset G$ be a subgroup isomorphic to $\kS_n$. Then for all $\sigma\in\kS_n$, there exists a unique $\lambda_\sigma\in\Gamma^n_0$ such that $(\lambda_\sigma,\sigma)\in H$.
\end{lemma}
Consider a transposition $(1,i)\in\kS_n$ with $1\leq i\leq n-1$. By the lemma, there exists $\lambda=(\lambda_1,\ldots,\lambda_n)\in\Gamma^n_0$ such that $(\lambda,(1,i))\in\Stab(F(0))$. This gives the relations:
$$
w_1-w_i=\lambda_1=-\lambda_i, \quad \lambda_2=\cdots=\widehat{\lambda_i}=\cdots=\lambda_n=0.
$$
Since $w_1$ and $w_i$ are in a fundamental domain of $\Gamma$, this forces $w_1=w_i$ and $\lambda=(0,\ldots,0)$. We obtain finally $w_1=\cdots=w_{n-1}=:w$ with $w_n=-(n-1)w$. Now consider the transposition $(1,n)$ and $(\lambda,(1,n))\in\Stab(F(0))$. This gives
$$
w-w_n=\lambda_1=-\lambda_n, \quad \lambda_2=\cdots=\lambda_{n-1}=0,
$$ 
so $nw\in\Gamma$. Let $a=-[w]$ be the class of $-w$ in $A=\IC^2/\Gamma$; then $na=~0$, so $a$ is an $n$-torsion point. Denote by $t_a\colon A\to A$ the translation morphism and by $t_a^{((n))}\colon\fiberA\to\fiberA$ the induced natural automorphism. Consider the automorphism $t_a^{((n))}\circ f$. One can choose the lift $T_a$ of $t_a^{((n))}$ to $(\IC^2)^n_0\to(\IC^2)^n_0$ to be the translation by $(-w,\ldots,-w,(n-1)w)$. Then $(T_a\circ F)(0)=0$.

\par\noindent{\it Step 3. Linearity.}\\ 
Let $(\gamma,\sigma)\in G$ and $(\gamma',\sigma'){\coloneqq}f_*(\gamma,\sigma)$. The $f_*$-equivariance of $F$ means
\begin{equation}\label{eq:F}
\forall z\in(\IC^2)^n_0, \quad F(\sigma z+\gamma)=\sigma'F(z)+\gamma'.
\end{equation}
Putting $z=0$, one gets $F(\gamma)=\gamma'$ (put in another way, $F(\Gamma^n_0)=\Gamma^n_0$). Moreover, observe that for any $\gamma\in\Gamma^n_0$, if $(\gamma'',\sigma''){\coloneqq}f_*(n!\gamma,\id)$ then $\sigma''=\id$. Indeed, setting $(\gamma',\sigma')=f_*(\gamma,\id)$ one computes
$$
(\gamma'',\sigma'')=f_*(n!\gamma,\id)=f_*\left((\gamma,\id)^{n!}\right)=(\gamma',\sigma')^{n!}=(\gamma'+\sigma'\gamma'+\ldots+(\sigma')^{n!-1}\gamma',\id).
$$
Putting $\sigma=\id$, in (\ref{eq:F}) one gets then for all $\gamma\in n!\Gamma^n_0$:
$$
\forall z\in(\IC^2)^n_0, \quad F(z+\gamma)=F(z)+F(\gamma).
$$
It follows that the partial derivatives of $F$ are $n!\Gamma^n_0$-periodic, hence constant by Liouville's theorem. Since $F(0)=0$, $F$ is linear.

\par\noindent{\it Step 4. Naturality.}\\
Let $(\gamma,\sigma)\in G$ and $(\gamma',\sigma'){\coloneqq}f_*(\gamma,\sigma)$. Since $F(0)=0$, the $f_*$-equivariance shows (putting $z=0$) that $\gamma=0$ if and only if $\gamma'=0$, so $f_*$ induces an isomorphism $\phi{\coloneqq}(f_*)_{|\kS_n}\colon\kS_n\xrightarrow{\sim}\kS_n$ and $F$ is $\phi$-equivariant:
$$
\forall \sigma\in\kS_n, \forall z\in(\IC^2)^n_0, \quad F(\sigma z)=\phi(\sigma) F(z).
$$

\par If $\phi$ is an inner automorphism (this is always the case when $n\neq 6$), by permuting the coordinates at the target, one can assume that
$\phi=\id$. Writting $F$ as a matrix $(a_{i,j})_{1\leq i,j\leq n}$ where the $a_{i,j}$ are $(2\times 2)$-matrices, the condition $\sigma^{-1}F\sigma=F$ gives easily (by evaluating at any transposition):
$$
a_{1,1}=\cdots=a_{n,n}=:\beta \text{ and } a_{i,j}=:\alpha \quad\forall i\neq j 
$$
for some matrices $\alpha,\beta$. Substituting $\delta{\coloneqq}\beta-\alpha$ (observe that $\alpha\neq\beta$ since $F$ is an isomorphism) and denoting $G{\coloneqq}(\alpha)_{1\leq i,j\leq n}$, one has $F=\diag(\delta,\ldots,\delta)+G$. Using the fact that for any $z\in(\IC^2)^n_0$, 
$$
G(z)=(\alpha(z_1+\ldots+z_n),\ldots,\alpha(z_1+\ldots+z_n))=(0,\ldots,0),
$$
one deduces that $F=\diag(\delta,\ldots,\delta)$ and $\delta\colon\IC^2\to\IC^2$ is an isomorphism. Since $F(\Gamma^n_0)=\Gamma^n_0$, one has $\delta(\Gamma)=\Gamma$, so $\delta$ defines a group automorphism of $A$ which induces $f$, so $f$ is natural.

\par If $n=6$, then $\phi$ could be an outer automorphism. Since two such automorphisms only differ by an inner automorphism, by the preceding discussion it is enough to consider the following automorphism, characterized by:
\begin{align*} 
\phi((1,2))&=(1,2)(3,4)(5,6),\quad \phi((2,3))=(1,4)(2,5)(3,6)\\
\phi((3,4))&=(1,3)(2,4)(5,6),\quad \phi((4,5))=(1,2)(3,6)(4,5)\\
\phi((5,6))&=(1,4)(2,3)(5,6)
\end{align*}
In this case, a quite long but elementary computation shows that all $(2\times 2)$-entries of $F$ are the same, so $F$ would not be an isomorphism: this case does not occur.
\end{proof}

\begin{corollary}\text{}
\begin{enumerate}
\item \label{cor:item1} The map $\Tors_n(A)\rtimes\GAut(A)\to\Aut(\genKm),\psi\mapsto \psi^{{\llbracket}n\rrbracket}$ is injective.
\item The kernel of the map $\Aut(\genKm)\to\cO(H^2(\genKm,\IZ)),f\mapsto f^*$ is isomorphic to $\Tors_n(A)\rtimes\IZ/2\IZ$.
\item If $A$ is a complex torus such that $\Pic(A)=\{0\}$, then every automorphism of $\genKm$ is natural.
\end{enumerate}
\end{corollary}

\begin{proof}\text{}
\par\noindent{(1) } Let $\psi\in\Tors_n(A)\rtimes\GAut(A)$ such that $\psi^{{\llbracket}n\rrbracket}=\id$. Decomposing $\psi=t_a\circ h$, one sees that necessarily $a=0$ (otherwise the subschemes supported at the origin would not be fixed), so $\psi=h$ and, as in the proof of Theorem~\ref{th:naturalKummer}, $\psi^{{\llbracket}n\rrbracket}$ is uniquely determined by the map $F=(h\times \cdots\times h)\colon(\IC^2)^n_0\to(\IC^2)^n_0$, so $F=\id$ and $h=\id$.

\par\noindent{(2) } Assume $f\in\Aut(\genKm)$ acts as the identity on $H^2(\genKm,\IZ)$. Since there is a decomposition $H^2(\genKm,\IZ)\cong H^2(A,\IZ)\oplus\frac{1}{2}\IZ[E_0]$ and $E_0$ is rigid, $f$ leaves $E_0$ globally invariant, so by Theorem~\ref{th:naturalKummer} it is a natural automorphism. Since the kernel of the map $\Aut(\genKm)\to\cO(H^2(\genKm,\IZ))$ is finite (see \S\ref{ss:IHS}), $f$ is of finite order. For any $a\in\Tors_n(A)$, the translation by $a$ on $A$ induces an automorphism $t_a^{{\llbracket}n\rrbracket}$ acting as the identity on $H^2(\genKm,\IZ)$: since the automorphism $t_a^{[n]}$ induced on $\hilbA$ is homotopic to the identity and the restriction map $H^2(\hilbA,\IC)\to H^2(\genKm,\IC)$ is surjective, $t_a^{{\llbracket}n\rrbracket}$ acts as the identity on $H^2(\genKm,\IC)$ (see Beauville \cite[\S 5, Proposition~9]{B}), but $H^2(\genKm,\IZ)$ has no torsion, so $t_a^{{\llbracket}n\rrbracket}$ acts as the identity on it. It remains to consider the case where $f$ comes from a group automorphism $h\in\GAut(A)$: $f=h^{{\llbracket}n\rrbracket}$. Let $d$ be the order of $f$. Then 
$(h^{{\llbracket}n\rrbracket})^d=(h^d)^{{\llbracket}n\rrbracket}=\id$ so by (\ref{cor:item1}), $h^d=\id$, one has that $h$ is determined by a $\IC$-linear isomorphism $H\colon\IC^2\to\IC^2$ such that $H^d=\id$. In a suitable $\IC$-basis $(v_1,v_2)$ of $\IC^2$, $H$ is given by a diagonal matrix $H=\diag(\lambda,\mu)$, so in the corresponding $\IR$-basis $(v_1,\overline{v_1},v_2,\overline{v_2})$ of $\Gamma\otimes_\IZ\IC$, $H$ is given by the diagonal matrix $\diag(\lambda,\overline{\lambda},\mu,\overline{\mu})$. Since $H^2(\genKm,\IC)\cong H^2(A,\IC)\oplus \IC[E_0]$, $h$ acts trivially on $H^2(A,\IC)\cong(\wedge^2\Hom_\IZ(\Gamma,\IC))\otimes_\IZ\IC$ but the action of $h$ is given by $\diag(\lambda\mu,\lambda\overline{\mu},\lambda\overline{\lambda},\overline{\lambda}\mu,\overline{\lambda}\overline{\mu},\mu\overline{\mu})$. The only possibility is $\lambda=\mu=\pm 1$, so $f=\pm\id$. 
\par\noindent{(3) } Recall that $\Pic(\genKm)\cong\Pic(A)\oplus\IZ\cL_0$ where $\cL_0^{\otimes 2}\cong\cO(-E_0)$. If $\Pic(A)=\{0\}$ then any automorphism of $\genKm$ leaves $E_0$ globally invariant, so is natural.
\end{proof}

\subsection{Topological Lefschetz numbers of the natural automorphisms}

The generalized Kummer variety $\genKm$ fits into the cartesian diagram
$$
\xymatrix{A\times\genKm\ar[r]^-\nu\ar[d]_p & \hilbA\ar[d]^S\\A\ar[r]_{n} & A}
$$
where $n\colon A \to A$ is the multiplication by $n$ and $\nu(a, \xi)= a + \xi$ is a Galois covering with Galois group $G {\coloneqq} \Tors_n(A)$. Here, $G$ acts via $g \cdot (a, \xi) = (a - g, g + \xi)$ on $A \times \genKm$.

Let $\psi\in\Aut(A)$, decomposed as $\psi=t_a\circ h$ with $a\in\Tors_n(A)$ and $h\in\GAut(A)$. 
If we let $\psi$ act on $A \times \genKm$ by $h\times\psi^{\llbracket n\rrbracket}\colon(a, \xi)\mapsto(h(a), \psi^{\llbracket n\rrbracket}(\xi))$ the cartesian diagram is equivariant with respect to $h\colon A \to A$ restricted to $G$. The topological Lefschetz number of $\psi^{\llbracket n\rrbracket}$ is by definition
$$
L(\psi^{\llbracket n\rrbracket}) = \sum_i (-1)^i \trace \left((\psi^{{\llbracket}n\rrbracket})^*|_{H^i(\Kummer n A, \IC)}\right)
$$
and one has the relation $L(h)\cdot L(\psi^{{\llbracket}n\rrbracket}) = L(h \times \psi^{{\llbracket}n\rrbracket})$.

By~\cite{NW}, there is a natural isomorphism
$$
	H^*(A \times \Kummer n A, \IC[2n]) \cong \left.\frac{\diff^n}{n! \diff\! t^n}\right|_{t = 0} \bigoplus_{\chi \in \dual G} 
	\Sym^* \left(\bigoplus_{\nu \ge 1} H^*(A, \IC[2]) \cdot t^{\nu \abs{\chi}}\right),
$$
where $\dual G$ is the dual group of $G$, $\abs\chi$ denotes the order of $\chi$ in $\dual G$ and $\Sym^*(\cdot)$ is the symmetric algebra. The action of \(h \times\psi^{{\llbracket}n\rrbracket}\) on the left hand side
is compatible with the action of $h$ on the right hand side. In particular, we may assume that $a = 0$ as $L(\psi^{{\llbracket}n\rrbracket})$
will be independent of $a$ since $t_a$ acts trivially on the cohomology. Note, in particular, that \(h\) also acts on \(\dual G\). The subgroup of invariant characters in
\(\dual G\) is denoted by \((\dual G)^h\).

\begin{lemma}
	Let \(h\colon H \to H\) be an even homomorphism between super vector spaces over a field \(K\). Then
	\[
		\str \Sym^* (h t) = \exp(-\str \log (1 - ht)) = \exp\left(\sum_{s \ge 1} \frac{\str h^s}{s} t^s\right),
	\]
	where \(\str\) denotes the super trace of an endomorphism between super vector spaces.
\end{lemma}

\begin{proof}
	By scalar extension to the algebraic closure, we may assume that \(K\) is algebraically closed.
	Furthermore, both sides are continuous in the Zariski topology on \(\End(H)\). Thus we may even assume that \(h\) is
	semi-simple.
	
	Assume that the claimed formula is true for \(h_1\colon H_1 \to H_1\) and \(h_2\colon H_2 \to H_2\). It follows that
	the formula holds for \(H = H_1 \oplus H_2\) and \(h = h_1 \oplus h_2\colon H \to H\) as
	\[
		\str \Sym^* ((h_1 \oplus h_2) t) = \str (\Sym^*(h_1 t) \otimes \Sym^*(h_2 t))
		= \str(\Sym^* (h_1 t)) \cdot \str (\Sym^* (h_2 t))
	\]
	and
	\[
		\begin{aligned}
		\exp(- \str \log (1 - (h_1 \oplus h_2) t)
		& = \exp(- \str (\log (1 - h_1 t) \oplus \log(1 - h_2 t)))
		\\
		& = \exp(- \str \log (1 - h_1 t) - \str \log (1 - h_2 t))
		\\
		& = \exp(- \str \log (1 - h_1 t)) \cdot \exp(- \str \log (1 - h_2 t)).
		\end{aligned}
	\]
	
	As we have assumed \(h\) to be semi-simple, it is diagonalisable. In view of the preceding, it
	is thus enough to prove the claimed formula in case \(H\) is a one-dimensional even or odd super vector space and
	\(h\) is multiplication by a scalar, which we call \(h\) again.
	
	In the even case, one has
	\[
		\str \Sym^* (h t) = \sum_{n = 0}^\infty h^n t^n = \frac 1 {1 - ht} = \exp(-\str \log(1 - ht));
	\]
	in the odd case, one has
	\[
		\str \Sym^* (h t) = 1 - ht = \exp(- \str \log(1 - ht)).
	\]
\end{proof}

One then has
\[
	\begin{aligned}
		L(h)\cdot L(\psi^{{\llbracket}n\rrbracket}) & =  L(h \times \psi^{{\llbracket}n\rrbracket})
		\\
		& =  \left.\frac{\diff^n}{n! \diff\! t^n}\right|_{t = 0}
		\sum_{\chi \in (\dual G)^h}
		\prod_{\nu \ge 1}
		L(\Sym^* (H^*(A, \IC[2]) \cdot t^{\nu \abs\chi}))
		\\
		& =\left.\frac{\diff^n}{n! \diff\! t^n}\right|_{t = 0}
		\sum_{\chi \in (\dual G)^h}
		\prod_{\nu \ge 1}
		\exp\left(\sum_{s \ge 1} \frac{\str((h^*)^s)} s t^{\nu \abs\chi s}\right),
	\end{aligned}
\]
where \(h^*\colon H^*(A, \IC[2]) \to H^*(A, \IC[2])\) is the induced graded homomorphism on
cohomology. Let \(\Psi = h^{-1}\colon H^{-1}(A, \IC[2]) \to H^{-1}(A, \IC[2])\). Since the cohomology of the torus is an exterior algebra over \(H^{-1}(A, \IC[2])\), one has
\(\str((h^*)^s) = \det (1 - \Psi^s)\). We have thus obtained the following formula:

\begin{proposition}\label{prop:lefschetz} Let $\psi\in\Aut(A)$ and $\Psi\coloneqq \psi^*\colon H^1(A,\IC)\to H^1(A,\IC)$. Then
\[
	\begin{aligned}
		L(\psi)L(\psi^{{\llbracket}n\rrbracket})
		& =  \left.\frac{\diff^n}{n! \diff\! t^n}\right|_{t = 0}
		\sum_{\chi \in (\dual G)^h}
		\prod_{\nu \ge 1}
		\exp\left(\sum_{s \ge 1} \frac{\det(1 - \Psi^s)} s t^{\nu \abs\chi s}\right).
	\end{aligned}
\]
\end{proposition}

\section{Construction of Enriques varieties}

To construct Enriques varieties of dimension $2n-2>2$ (so $n>2$), a first guess could be to take the Hilbert scheme
of $n-1$ points on a K3 surface $S$ (since $\dim S^{[n-1]}=2n-2$) and look for an automorphism of order $n>2$ acting fixed point free. Apart from Beauville's involution on the Hilbert scheme of two points on a generic quartic surface (that has fixed points, as is easily seen), the only known automorphisms are the natural automorphisms, given by the natural map $\Aut(S)\to\Aut(S^{[n-1]})$, $\psi\mapsto \psi^{[n-1]}$  \cite{Boissiere,BS}. Since an automorphism $\psi$ of $S$ is (non)-symplectic if and only if $\psi^{[n-1]}$ is, one could start from a purely non-symplectic automorphism of order $n>2$ on $S$. But the holomorphic Lefschetz number of such an automorphism is never zero, so $\psi$ has fixed points. Then any point of $S^{[n-1]}$ whose scheme structure is supported at such a fixed point and whose ideal is monomial (in a coordinate system for which the action of $\psi$ is linearized and diagonal) is a fixed point for $\psi^{[n-1]}$ \cite[Proposition 7]{Boissiere} thus no higher dimensional Enriques variety can be constructed like this. 

\subsection{Fixed points of natural automorphisms on generalized Kummer varieties}

Let $A=\IC^2/\Gamma$ be a complex torus, $n\geq 3$ and $\psi=t_a\circ h$ as before, where $a\in\Tors_n(A)$. We are looking for automorphisms $\psi$ of finite order such that $\psi^{{\llbracket}n\rrbracket}$ has no fixed point on $\genKm$. We can already exclude two trivial cases:
\begin{itemize}
\item If $\psi=h\in\GAut(A)$, then $\psi(0)=0$ so any subscheme of $\genKm$ supported at the origin, whose ideal is monomial (in a coordinate system for which the action of $\psi$ is linearized and diagonal) is fixed by $\psi^{{\llbracket}n\rrbracket}$.
\item If $\psi=t_a$, then $\psi$ has no fixed points on $A$ but for any $n\geq 3$, the subschemes consisting of $n$ points $\{x,x+a,\ldots,x+(n-1)a\}$ for appropriate values of $x$ are fixed by $\psi^{{\llbracket}n\rrbracket}$. For odd $n$, take any point $x\in\Tors_n(A)$ ; for even $n$, let $\alpha\in A$ such that $2\alpha=a$ and take $x=\alpha+y$ with $y\in\Tors_n(A)$. 
\end{itemize}
Moreover, if $\psi\in\Tors_n(A)\rtimes\GAut(A)$ is of finite order $n\geq 2$ on a simple torus $A$ and if $\Fix(\psi)=\emptyset$, then $\psi^{{\llbracket}n\rrbracket}$ has fixed points on $\genKm$: on a simple torus, any biholomorphic map without fixed points is a translation \cite[\S13.1]{BL}, so $\psi^{{\llbracket}n\rrbracket}$ has fixed points on $\genKm$ by the preceding observation.
But observe that even though a given automorphism $\psi$ of order $n$ on a (non-simple) complex torus has no fixed points, the corresponding automorphism $\psi^{{\llbracket}n\rrbracket}$ on $\genKm$ can still have fixed points, for instance orbits $\{x,\psi(x),\ldots,\psi^{n-1}(x)\}$ whose sum is zero.

Conversely, if an automorphism $\psi\in\Aut(A)$ has fixed points, it is not always easy to see if $\psi^{{\llbracket}n\rrbracket}$ has fixed points or not. In some cases, this question can be solved by computing its topological Lefschetz number by Proposition~\ref{prop:lefschetz}. Consider for example an automorphism $\psi$ of order $5$ on a two-dimensional torus $A$ and its action on $\Kummer 5 A$ and denote as usual by $h$ its linear part. By the classification
of these automorphisms on two-dimensional complex tori \cite{BL}, \(\Psi\) (see notation in Proposition~\ref{prop:lefschetz}) is given by a diagonal matrix whose entries
are the primitive \(5\)-roots of unity, \(\xi, \dotsc, \xi^4\). It follows that \(L(\psi) = 5\).
The group \(\dual G\) has one element of order \(1\) which is,
of course, a fixed point under \(h\) and \(624\) elements of order \(5\), of which \(4\) elements are also fixed points of \(h\). (There
are in total \(5\) fixed points of \(h\) on \(A\) and all of them are \(5\)-torsion points.) It follows that
\[
	\begin{aligned}
		L(\Kummer 5 A)
		& = \left.\left(\frac 1 5 \frac{\diff^5}{5! \diff\! t^5}
		+ \frac 4 5 \frac{\diff}{\diff\! t}\right)\right|_{t = 0}
		\prod_{\nu \ge 1}
		\exp\left(\sum_{s \ge 1} \frac{\det(1 - F^s)} s t^{\nu s}\right).
	\end{aligned}
\]
This calculation can be done explicitely. Using \(\det (1 - F^s) = 5\) in case \(s\) is not divisible by \(5\) and
\(\det (1 - F^s) = 0\) otherwise, one has
\[
	\begin{aligned}
		\exp\left(\sum_{s \ge 1} \frac{\det(1 - F^s)} s t^{\nu s}\right)
		& = \exp\left(\sum_{s \ge 1} \frac{5 t^{\nu s}} s - \sum_{s \ge 1} \frac {5 t^{5 \nu s}}{5 \nu s}\right)
		\\
		& = \exp(- 5 \log(1 - t^\nu)) \exp(\log(1 - t^{5 \nu}))
		= \frac{1 - t^{5 \nu}}{(1 - t^\nu)^5},
	\end{aligned}
\]
so that
\[
	\begin{aligned}
		L(\Kummer 5 A)
		& = \left.\left(\frac 1 5 \frac{\diff^5}{5! \diff\! t^5}
		+ \frac 4 5 \frac{\diff}{\diff\! t}\right)\right|_{t = 0}
		\prod_{\nu \ge 1}
		\frac{1 - t^{5\nu}}{(1 - t^\nu)^5}
		\\
		& = \left.\left(\frac 1 5 \frac{\diff^5}{5! \diff\! t^5}
		+ \frac 4 5 \frac{\diff}{\diff\! t}\right)\right|_{t = 0}
		\left(
			1 + 5 t + 20 t^2 + 65 t^3 + 190 t^4 + 505 t^5 + \mathrm O(t^6)
		\right)
		\\
		& = \frac 1 5 \cdot 505 + \frac 4 5 \cdot 5 = 105.
	\end{aligned}
\]

In particular, no automorphism of order \(5\) on a two-dimensional complex torus induces a fixed point free natural automorphism of a generalised Kummer variety of dimension \(8\).

\subsection{Construction of Enriques varieties of dimension $4$ and $6$}

\begin{proposition}\label{prop:construction}
Let $E$ be an elliptic curve admitting an automorphism of order $n\in\{3,4\}$, $\zeta_n$ a primitive $n$-th root of  unity and $A{\coloneqq}E\times E$. Let $h\in\GAut(A)$ be given by $h=\left(\begin{matrix}\zeta_n & 0 \\ 0 & 1\end{matrix}\right)$, $a_1,a_2\in E$(non zero) points of order $n$ in $E$, $a{\coloneqq}(a_1,a_2)$ and $\psi{\coloneqq}t_a\circ h$. Then for an appropriate choice of $a_1$, $\langle \psi^{{\llbracket}n\rrbracket}\rangle$ acts freely on $\genKm$ and $\genKm/\langle \psi^{{\llbracket}n\rrbracket}\rangle$ is an Enriques variety of dimension $2n-2$. 
\end{proposition}

\begin{proof} One has $E\cong\IC/(\IZ\oplus\zeta_n\IZ)$. Let $z=(x,y)\in A$. Then ${\psi(z)=(\zeta_nx+ a_1,y+a_2)}$ is an automorphism of order $n$ on $A$. Looking at the second coordinate, one sees that $\psi$ has no fixed point on $A$ and that there is no orbit of length strictly smaller than $n$ (since $a_2$ has order $n$). We study all possible fixed points for $\psi^{{\llbracket}n\rrbracket}$ (and its iterates) on $\genKm$.
\par\noindent{{(1)}} $n=3$. The only possibility is a fixed point on $\Kummer 3 A$ whose support is an orbit $\{z,\psi(z),\psi^2(z)\}$ such that the sum of the points is zero. Looking at the first coordinate, one gets the condition:
$(2+\zeta_3)a_1=0$. Taking for example $a_1=1/3$, $(2+\zeta_3)/3\neq 0$ in $E$, so $\psi^{{\llbracket}3\rrbracket}$ has no fixed points.

\par\noindent{{(2)}} $n=4$. The only possible fixed points for $\psi^{{\llbracket}4\rrbracket}$ are supported on length four orbits $\{z,\psi(z),\psi^2(z),\psi^3(z)\}$ whose sum is zero. The condition on the first coordinate is $2(1+\ii)a_1=0$, so if $a_1=1/4$ then $\psi^{{\llbracket}4\rrbracket}$ has no fixed points. Similarly, $(\psi^{{\llbracket}4\rrbracket})^2=(\psi^2)^{{\llbracket}4\rrbracket}$ can have fixed points supported on $\{z,\psi^2(z)\}\cup\{w,\psi^2(w)\}$ (for $z=w$, this means that the fixed point on $\Kummer 4 A$ has a non-reduced scheme structure). One gets for the first coordinate the same condition $2(1+\ii)a_1=0$, so for $a_1=1/4$ the group $\langle \psi^{{\llbracket}4\rrbracket}\rangle$ acts freely on $\Kummer 4 A$.
\end{proof}

\begin{remark}\label{rem:intermediate}
For $n=6$, in order to avoid fixed points supported on orbits of length $6$ on $\Kummer 6 A$, the same construction gives the condition $6\zeta_6a_1=0$ on the first coordinate. Unfortunately, every order-$6$ point $a_1\in E$ satisfies this equation, so the similarly defined non-symplectic automorphism on $\Kummer 6 A$ does have fixed points. This means that the quotient $\Kummer 6 A /\langle f^{{\llbracket}6\rrbracket}\rangle$ is not an Enriques variety. But consider again the automorphism of order three $f(z)=(\zeta_3x+a_1,y+a_2)$, which also acts on $\Kummer 6 A$. Its only possible fixed points are supported on union of orbits $\{z,f(z),f^2(z)\}\cup\{w,f(w),f^2(w)\}$ (where $z=w$ is possible). Again, since the sum of the points is zero, for the first coordinate the condition is $2(2+\zeta_3)a_1=0$, so for example if $a_1=1/3$, then $f^{{\llbracket}6\rrbracket}$ has no fixed point on $\Kummer 6 A$: the quotient $\Kummer 6 A/\langle f^{{\llbracket}6\rrbracket}\rangle$ is then a weak Enriques variety.
\end{remark}

\subsection{A counter-example to Theorem \ref{prop:classification}(\ref{prop:itemirredodd}) for even $d$}\label{ss:counterexample}

We construct a ten-dimensional Enriques variety that is not the quotient of an irreducible holomorphic symplectic manifold by an automorphism acting freely.

Take $W{\coloneqq}\Kummer 3 A$ with an automorphism $f^{{\llbracket}3\rrbracket}$ of order three, constructed in Proposition \ref{prop:construction}. Let $[x_0:\ldots:x_7:y_0:\ldots:y_7]$ be coordinates in $\IP_{13}$ and consider the intersection of seven quadrics of equations $Q_j(x)-Q'_j(y)=0$. For a generic choice, this intersection is complete and smooth and is a six-dimensional Calabi--Yau variety $V$. Consider the involution $\iota$ on $\IP_{13}$ such that $\iota(x_j)=-x_j$ and $\iota(y_j)=y_j$ for all $j$. The fixed locus of $\iota$ in $\IP_{13}$ is the disjoint union of two six-dimensional projective spaces, so for a generic choice of the quadrics it does not intersect $V$, hence $\iota$ acts freely on $V$. Consider now $X=W\times V$ with the automorphism $f^{{\llbracket}3\rrbracket}\times\iota$. It is easy to see that the order six group generated by $f^{{\llbracket}3\rrbracket}\times\iota$ acts freely on $X$ and that $Y{{\coloneqq}}X/\langle f^{{\llbracket}3\rrbracket}\times~\iota\rangle$ is an Enriques variety of dimension $10$ that can not be the unramified cyclic quotient of an irreducible holomorphic symplectic manifold. Observe however that $Y\cong \left(W/\langle f^{{\llbracket}3\rrbracket}\rangle\right)\times \left(V/\langle\iota\rangle\right)$ so in fact this Enriques variety is reducible.

\bibliographystyle{amsplain}
\bibliography{BiblioEnriquesKummer}

\end{document}